\documentclass{amsart}
\usepackage{latexsym}
\usepackage{amssymb}

\begin{document}
\title{REAL VALUED DIMENSION FOR MORSE-SARD THEOREM}
\author{Azat Ainouline}
\address{}
\email{azat.ainouline@3web.net}
\begin{abstract}
A dimension allowing in particular to state necessary and sufficient conditions of the Morse-Sard theorem for real valued function is introduced.
\end{abstract}
\maketitle

\newcommand{\real}{\mathbb R}
\newcommand{\meas}{{\rm meas}}
\newtheorem{theorem}{Theorem}
\newtheorem{corollary}{Corollary}
\newtheorem{lemma}{Lemma}[section]
\newcommand{\diam}{{\rm diam}}
\newcommand{\rank}{\rm rank}
\newtheorem{definition}{Definition}[section]
\numberwithin{equation}{section}

\section{Introduction.}
The aim of this paper is to find a real valued dimension that is a good extention of the standard integer valued dimension on  subsets of Euclidean spaces, in a sense that it would have  more "good" properties, as for instance can be used to state necessary and sufficient conditions for the Morse-Sard Theorem.

A number of real valued dimensions have been intensively studied  in geometric measure theory (see for example Mattila \cite{Mattila}), and some of the dimensions have desired properties such as : $\dim R^n=n, ~A\subseteq B\Longrightarrow \dim A \leqslant \dim B,~\dim\bigcup_{i\in \mathbb N}A_i=\sup_{i\in \mathbb N}A_i$~~etc. One of the central dimensions - the Hausdorff dimension - was used by many authors for sharpening the necessary part of the Morse-Sard theorem , for example:\newline
{\bf Theorem (Federer \cite{Federer})}\newline\it
Let $F:\real^n\rightarrow \real^m $~~be a $C^k$~~function,~~$k\in \mathbb N,~~ C_p(F)=\{x\in \real^n | rank (DF)\leqslant p\}$,~~then the Hausdorff measure of dimension $p+\frac{n-p}{k}$~~of $F(C_p(F))$~~is zero.\\ 
\rm

 But it was shown that the Hausdorff dimension has lack of the property that $\dim_H (A\times B)=\dim_H A+\dim_H B$~~(see Federer \cite{Federer} 2.10.29), and as we will see later, it does not give "right" dimension for set of critical values in the Morse-Sard theorem in sense that one cannot state that: \newline
\it if Hausdorff measure of dimension ~~$p+\frac{n-p}{k}$~~of a set $A$~~is zero then $A\subseteq F(C_p(F))$~~for some ~~$C^k$~~function $F:\real^n\rightarrow \real^m$\rm(see theorems  \ref{p14},\ref{p9},\ref{p15} here).\\

 In this article the author introduces a promising candidate for the "right" real valued dimension that can fix at least that weakness of the Hausdorff dimension, regarding to a case of real valued function in the Morse-Sard theorem (theorems  \ref{p9},\ref{p15}). 

\begin{definition}\label{D^k}\rm
For $m,n\in\mathbb N,~k\in\real$~~a  function $\psi:B\subseteq\real^m\rightarrow \real^n$~~is {\bf {\boldmath $D^k$}-function}~ if 
$\exists K>0~:~~\forall b,b^\prime\in B \quad|\psi(b)-\psi(b^\prime)|^k\leqslant K|b-b^\prime|.$
\end{definition}

\begin{definition}\label{R-dimension}\rm
For $n\in\mathbb N$~~  {\bf Outer dimension} of the set $A\subseteq \real^n$ is a number $Dm(A)=\inf\{k\in\real;~\exists\psi:B\subseteq\real^1\stackrel{D^k}{\rightarrow}\real^n,~k\geqslant 0,~A\subseteq\psi(B)\}.$
\end{definition}

\noindent

\begin{theorem}\label{p11}
Let $n\in \mathbb N$  then:
\begin{enumerate}
\item~~$Dm(\emptyset)=0$
\item~~$Dm(A)\leqslant Dm(B)$~~whenever  $A\subseteq B\subseteq\real^n$
\item~~$Dm(\bigcup_{i\in \mathbb N} A_i)~=~\sup_{i\in \mathbb N} Dm (A_i)$~~whenever $A_1,A_2,....\subseteq \real^n$.
\end{enumerate}
\end{theorem}

\begin{theorem}\label{p2}
Every $m$-dimensional $C^1$ manifold $A\subseteq\real^n,~(m,n\in\mathbb N,~m\leqslant n)$~~has $Dm$-dimension $m$.
\end{theorem}

\begin{theorem}\label{p5}
$ Dm(\psi(B))\leqslant k\cdot Dm(B)$~for any $D^k$-function $\psi:B\subseteq\real^m\rightarrow\real^n,~k\in\real$.
\end{theorem}

\begin{definition}\rm
For $A\subseteq\real^n$ {\bf Hausdorff dimension} of $A$ is a number 
$\dim_H A=\inf\left\{k;~\lim\limits_{\sigma\to 0}\inf\{\sum^\infty_{i=1}(\diam(B_i))^k;~B_i\subseteq\real^n,~\diam 
(B_i)<\sigma,~A\subseteq\bigcup^\infty_{i=1}B_i\}<\infty\right\}.$
\end{definition}

\begin{theorem}\label{p14}
\begin{enumerate}
\item~~$\dim_H\leqslant Dm$
\item~~$\exists A\subseteq\real^1$- compact, such that $\dim_H A< Dm (A)$
\end{enumerate} 
\rm\end{theorem}

\begin{theorem}\label{p4}
If $A, B\subseteq\real^n,~(n\in\mathbb N),$~then $ Dm(A\times B)\leqslant Dm(A)+Dm(B).$
\end{theorem}

\begin{corollary}[from Theorem \ref{p4} and Marstrand's Theorem \cite{Marstrand}]\label{c4}
If $A, B\subseteq\real^n,~(n\in\mathbb N),$~then $\dim_H A + \dim_H B\leqslant Dm(A\times B)\leqslant Dm(A)+Dm(B).$
\end{corollary}

\begin{theorem}\label{p1}
If $Dm(A)\not=Dm(B)$,~~ then~~sets $A$~and~$B$~ are not diffeomorphic.
\end{theorem}

\begin{definition}\label{Holder}\rm
For $n,m\in \mathbb N,k\in [0,\infty)$~~ a function $f:\real^m\rightarrow\real^n$~~ satisfies a {\bf $k$-H\"{o}lder condition}~~(or $f\in H^k $)~~if for every compact neighborhood $U$ there exists $M>0$~ such that 
$|f(x)-f(y)|\leqslant M\cdot|x-y|^k \text{~~for all~~} x,y\in U.$
\end{definition}
\noindent
The next theorem is another definition of the $Dm$~~dimension.
\begin{theorem}\label{Dm}
Let $A\subseteq \real^n,~~n\in \mathbb N$~~then:~~\newline
$Dm(A)=\inf\{1/k;~\exists\psi:B\subseteq\real^1\stackrel{H^k,onto}{\longrightarrow} A \}=\inf\{k;~\exists\psi:B\subseteq\real^1\stackrel{H^{1/k},onto}{\longrightarrow} A \}.$
\end{theorem}

\begin{definition}\rm
For $k=0,1,2...,~~\lambda \in[0,1)$~~a function~~$ f:\real^m\rightarrow \real^n$~~is a  $C^{k.\lambda}$-function (or {\boldmath $f\in C^{k.\lambda}$}),~~if $f\in C^k$ and every $k^{th}$ partial derivative of $f$  is a $H^\lambda$-function. 
\end{definition}

Next three theorems show relationship among dimensions, $C^k$ functions, critical points and critical values.

\begin{theorem}\label{p9}
Let $F:\real^n\stackrel{C^{k.\lambda}}{\rightarrow}\real^m,~~k\in \mathbb N,~~\lambda \in [0,1), ~~p<\min\{m,n\}$~ an integer.\quad Then 
$Dm(F(C_p(F)))\leqslant min\left\{p+\frac{n-p}{k+\lambda},m\right\}$.
\end{theorem}

\begin{theorem}\label{p15}
For any $A\subseteq \real^1,~(Dm(A)<1)$,~natural $m$ and real $k~(1\leqslant k <\frac{m}{Dm(A)}),~A$ is subset of the critical values set of some $C^{[k].k-[k]}$ function $f:\real^m\rightarrow\real^1$,~where $[k]$ is integer part of $k$.
\end{theorem}

Note: In the special case when $Dm(A) =1,$~the existence of appropriate function for Theorem \ref{p15} is shown in  \cite{Bates 1}.

\begin{theorem}\label{p7}
$Dm(f(\sum_{k+\lambda} f))\leqslant\frac{Dm(B)}{k+\lambda}$~for any $C^{k.\lambda}$~map~$f:B\subseteq\real^m\rightarrow\real^n~~(m,n\in\mathbb N,~k+\lambda \not= 0).$~~where $\sum\nolimits_{k+\lambda} f=\{x\in B;~\text{every partial derivative of}~ f~~\text{of order}~< k+\lambda~\text{vanishes at}~x\}$
\end{theorem}

Next we define  {\bf Outer {\boldmath$Dm$}-measure} $m_a,~ a\in [0,+\infty)$.

\begin{definition}\label{m_a}\rm Let $A\subseteq \real^n,~n\in \mathbb N ,~a\in [0,+\infty)$.~ Then if $a\geqslant Dm(A)$, \quad {\boldmath$m_a(A)$}$=\inf \{meas (B);~\exists f:B\subseteq \real \rightarrow \real^n,~ A\subseteq f(B),~ \forall b,b^\prime\in B \quad|\psi(b)-\psi(b^\prime)|^a\leqslant |b-b^\prime|\},$
~where $meas(B)$ is the Outer Lebesgue measure of $B$; otherwise $m_a(A)=\infty$.
\end{definition}

\noindent
To show correctness of the definition \ref{m_a} with regard to Dm-dimensions we will prove the following two theorems:
\begin{theorem}\label{p22}
Let $n\in \mathbb N,a\in [0,+\infty)$~~then:
\begin{enumerate}
\item $m_a(\emptyset)=0$
\item $m_a(A)\leqslant m_a(B)$ \quad whenever $A\subseteq B\subseteq \real^n$
\item $m_a(\bigcup_{i\in \mathbb N}A_i)\leqslant \sum_{i\in \mathbb N}m_a(A_i)$ \quad whenever $A_1,A_2,....\subseteq \real^n$.
\end{enumerate}
\end{theorem}

\begin{theorem}\label{p21}
Let $A\subseteq \real^n~~(n\in \mathbb N)$,~~then for any $a>Dm(A),   ~~m_a(A)=0$.
\end{theorem}

\section{Proofs.}

\noindent
{\bf Proof of Theorem \ref{p5}.}
Let $t\in\real$ be such that $Dm(B)<t$, then there exists a $D^t$ function $f:A\subseteq\real^1\rightarrow \real^m$, that $B\subseteq F(A)$. By the "Composition property" \cite{Ainouline2} ~~ $\psi\circ f:A\subseteq\real^1\rightarrow\real^n$ is a $D^{tk}$ function and also $\psi(f(A))\supseteq\psi(B)$. It means that $Dm(\psi(B))\leqslant k\cdot t,$ and because it is true for every $t>Dm(B)$, then $Dm(\psi(B))\leqslant k\cdot Dm(B).$
\hspace*{\fill}$\Box$ \\

\begin{lemma}\label{lemma 3.2}
Let $A=\bigcup_{i\in\mathbb N}A_i\subseteq\real^n,~n\in\mathbb N,~k\geqslant 0,~k\in\real$ and $\{f_i:B_i\subseteq[a_i,b_i]\stackrel{D^k}{\rightarrow}\real^n,~i\in\mathbb N\}$ be such that $A_i\subseteq f_i(B_i)$. Then there exists $f:B\subseteq \real^1\stackrel{D^k}{\rightarrow}\real^n$ with $A\subseteq f(B).$
\rm\end{lemma}

\noindent
{\bf Proof.}
By the definition of $D^k$-function there exist real positive numbers $K_i~(i\in \mathbb N),$ such that for each $i \quad |f_i(x)-f_i(y)|^k\leqslant K_i|x-y|. $ Also from  the "Extension on closure" property \cite{Ainouline2}~~every $f_i(B_i)$ is bounded set in $\real^n.$  Let us take sets $\{[x_i,y_i]\subseteq\real^1,~i\in\mathbb N\}$ and $\{l_i:[x_i,y_i]\stackrel{onto}{\rightarrow}[a_i,b_i],~i\in\mathbb N\}$ such that:
\begin{enumerate}
\item $x_1 =0$
\item $\forall i\in\mathbb N$
  \begin{itemize}
  \item $x_i<y_i<x_{i+1}$
  \item $y_i-x_i\geqslant K_i(b_i-a_i)$
  \item $x_{i+1}-y_i>\sup\limits_{j\leqslant i,a\in[a_j,b_j],~b\in[a_{i+1},b_{i+1}]} |f_j(a)-f_{i+1}(b)|^k$ 
  \end{itemize}
\item $\forall i\in\mathbb N,~\forall x\in[x_i,y_i] \quad l_i(x)=\frac{x-x_i}{y_i-x_i}(b_i-a_i)+a_i.$
\end{enumerate}

Now by defining $B=\bigcup l_i^{-1}(B_i)$ and $f:B\subseteq\real^1\rightarrow\real^n$ such that $f=f_i\circ l_i\upharpoonright (B\cap [x_i,y_i])$ we obtain the required $D^k$ function $f$ such that\newline
$|f(x)-f(y)|^k\leqslant 1\cdot|x-y| \quad \forall x,y\in B,\quad A\subseteq f(B).$
\hspace*{\fill}$\Box$ \\

\noindent
{\bf Proof of Theorem \ref{p11}.}
The parts 1),2) of this theorem are obvious. Let's prove the part 3).
$Dm(\bigcup_{i\in \mathbb N}A_i)\geqslant\sup_{i\in\mathbb N}Dm (A_i)$ follows from the part 2). Then suffice it to prove that $Dm(\bigcup_{i\in \mathbb N}A_i)\leqslant\sup_{i\in\mathbb N}Dm (A_i)$. Let $k\in\real$ be such that $\sup_{i\in\mathbb N}Dm (A_i)< k$, then for each $i$ there exists a $D^k$ function $g_i: M_i\subseteq \real^1\rightarrow\real^n$  such that  $A_i\subseteq g_i(M_i).$

Let us define a denumerable set of $D^k$ functions:
$F=\{f_{ij}=g_i\upharpoonright B_{ij};~~ B_{ij}=[j,j+1]\cap M_i,~j\in\mathbb Z,~i\in\mathbb N\}.$
Then we have a denumerable set $P=\{f_{ij}(B_{ij});~i\in\mathbb N,~j\in\mathbb Z\}$ such that $\bigcup_{i\in\mathbb N}A_i=\bigcup_{i\in\mathbb N,j\in\mathbb Z}f_{ij}(B_{ij})$. And the sets $F,P$ satisfy the condition of Lemma \ref{lemma 3.2}. So that there exists a function $f:B\subseteq\real^1\stackrel{D^k}{\rightarrow}\real^n$ with $\bigcup_{i\in\mathbb N} A_i\subseteq f(B)$. But it means that $Dm(\bigcup_{i\in\mathbb N}A_i)\leqslant k$. And because $Dm(\bigcup_{i\in \mathbb N}A_i)\leqslant k$~~for all $k$~~such that $\sup_{i\in\mathbb N}Dm (A_i)< k$,~~ it follows that $Dm (\bigcup_{i\in\mathbb N}A_i)\leqslant\sup_{i\in\mathbb N}Dm (A_i).$
\hspace*{\fill}$\Box$ \\

\noindent
{\bf Proof of Theorem \ref{Dm}.}
The proof based on the facts that every $D^k,~~k>0$~~is a $H^{1/k}$~~function, and on the other hand if $f\in H^{1/k}$~~then $\forall [a_i,b_i]\subseteq \real$~~a function $f_i=f\upharpoonright [a_i,b_i]$~~is a $D^k$~~function. Then by applying the lemma \ref{lemma 3.2}~~we finish the proof of Theorem \ref{Dm}.  
\hspace*{\fill}$\Box$\\

\noindent
{\bf Proof of Theorem \ref{p4}.} 
We  set ~$K^n_0=\{Q^n_{i_0},~i_0\in\mathbb N\}$,~where $Q^n_{i_0}$ is a closed cube in $\mathbb R^n$ with side length 1, and every coordinate of any vertex of $Q^n_{i_0}$~~is an integer. 
In general, having constructed the cubes of $K^n_{s-1}$,  divide each $Q^n_{i_0,i_1,i_2,....,i_{s-1}}\in K^n_{s-1}$ into $2^n$ closed cubes of side $\frac{1}{2^s}$,  and let $K^{n}_{s}$ be the set of all these cubes. More precisely we will write for $s\in \mathbb N$:
\begin{equation}
\label{K^n_s}
K^{n}_{s}=\{ Q^n_{i_0,i_1,i_2,...,i_{s-1},i_s}~~;~~ Q^n_{i_0,i_1,i_2,...,i_{s-1},i_s} \subseteq Q^n_{i_0,i_1,i_2,...i_{s-1}}\in K^n_{s-1}, 1\leqslant i_s\leqslant 2^n \}.
\end{equation}
And we apply here:
\begin{lemma}[\cite{Ainouline3}]\label{n,p,k}
Let $n,p\in \mathbb N,~p\leqslant n,~k\in\real,k\geqslant 1.$  Then there exists a continuous space-filling function
$\pi^n_{k,p}=(\pi_1,\pi_2):[0,1]\stackrel{onto}{\rightarrow}[0,1]^n$
such that $\pi_1:[0,1]\stackrel{D^{\frac{pk+n-p}{k}}}{\longrightarrow}[0,1]^p$  \quad and  \quad$\pi_2:[0,1]\stackrel{D^{pk+n-p}}{\longrightarrow}[0,1]^{n-p}$
are the component functions of $\pi^n_{k,p}.$
\end{lemma}
Now let $k,t$ be real numbers such that $Dm(A)<k,~Dm(B)<t$, then there exists a $D^k$ function $\psi^*:A^*\subseteq\real^1\rightarrow\real^n$ with $A\subseteq\psi^*(A^*)$ and a $D^t$ function $\psi^{**}:B^*\subseteq\real^1\rightarrow\real^n$ with $B\subseteq\psi^{**}(B^*).$
Let us consider the following:
\begin{itemize}
\item[{a)}] A function ~~$\psi:A^*\times B^*\subseteq\real^2\rightarrow\real^{2n}$ such that $\psi(x,y)=(\psi^*(x),\psi^{**}(y))\in \real^{2n} \quad\forall (x,y)\in A^*\times B^*$.
\item[{b)}] The set $K^2_0$ (see (\ref{K^n_s})).
\item[{c)}] A family of functions $\psi_\delta=\psi\upharpoonright\delta\cap(A^*\times B^*),~~\delta\in K^2_0$.
\item[{d)}] The function $\pi^{2}_{k/t,1}=(\pi_1,\pi_2):[0,1]\stackrel{onto}{\rightarrow}[0,1]^2$ (see Lemma \ref{n,p,k}).
\item[{e)}] A family of the identity maps $I_\delta:[0,1]^2\stackrel{onto}{\rightarrow}\delta~~(\delta\in K^2_0)$ such that $I_\delta=(I^*_\delta,I^{**}_\delta)$ and $I^*_\delta, I^{**}_\delta:[0,1]\stackrel{onto}{\rightarrow} \delta^*,\delta^{**} \in K^1_0$ are identity maps.
\item[{f)}] A family of sets $C_\delta=(\pi^{2}_{k/t,1})^{-1} (I^{-1}_\delta(\delta\cap(A^*\times B^*))),~~~\delta\in K^2_0.$
\item[{g)}] A family of functions $f_\delta:C_\delta\subseteq[0,1]\rightarrow\real^2~~~(\delta\in K^2_0),$ such that $\forall\delta\in K^2_0,~x\in C_\delta \quad f_\delta(x)=\psi_\delta(I_\delta(\pi^{2}_{k/t,1}(x))).$ 
\end{itemize}

Let us prove that \newline
1)~~$f_\delta$ is $D^{k+t}$ function for each $\delta \in K^2_0$\newline
2)~~$A\times B\subseteq\bigcup_{\delta\in K^2_0}f_\delta(C_\delta).$\\

\begin{itemize}
\item[{1)}] $\forall \delta\in K^2_0$~by Lemma \ref{n,p,k}~~and g) \quad $f_\delta=\psi_\delta\circ I_\delta\circ \pi^2_{k/t,1}=\psi_\delta\circ I_\delta\circ(\pi_1,\pi_2)=\psi_\delta\circ(I^*_\delta\circ\pi_1,I^{**}_\delta\circ\pi_2)= (\psi^*\circ I^*_\delta\circ\pi_1,\psi^{**}_\delta\circ I^*_\delta\circ\pi_2),$ where $\psi^*_\delta,~\psi^{**}_\delta$ are the component functions of $\psi_\delta$.\\
Recalling that $\psi^*\in D^k,~~\psi^{**}\in D^t$, we can see that $\psi^*_\delta\in D^k,~ \psi^{**}_\delta \in D^t, \quad \pi_1\in D^{\frac{k/t+2-1}{k/t}},~\pi_2\in D^{k/t+2-1}$ by Lemma \ref{n,p,k}, and $I^*_\delta,I^{**}_\delta\in D^1$ as identity maps.\\
Thus by "Composition property" \cite{Ainouline2}  $\psi^*\circ I^*_\delta\circ\pi_1\in D^{k+t}$ and $\psi^{**}\circ I^*_\delta\circ\pi_2\in D^{k+t}$. Consequently $f_\delta\in D^{k+t}$ as both its component functions are $D^{k+t}.$
\item[{2)}] Recall that $A\subseteq\psi^*(A^*),~B\subseteq\psi^{**}(B^*)$, then
  \begin{itemize}
  \item  $A\times B\subseteq rang~\psi$ \quad by {\bf a)};
  \item  $rang~\psi\subseteq\bigcup_{\delta\in K^2_0} rang~\psi_\delta$ \quad by {\bf c)}.
  \end{itemize}
Also $I_\delta:[0,1]^2\rightarrow\delta$ is "onto" as identity map and $\pi^2_{k/t,1}$ is "onto" by Lemma \ref{n,p,k}. So that by {\bf g)} $A\times B\subseteq\bigcup_{\delta\in K^2_0} rang f_\delta=\bigcup_{\delta\in K^2_0} f_\delta(C_\delta).$
\end{itemize}
Now from 1) and 2) we can see that the set $A\times B$ and the denumerable family of functions $f_\delta,~\delta\in K^2_0$ satisfy the conditions of Lemma \ref{lemma 3.2}, so that there exists $f:M\subseteq\real^1\stackrel{D^{k+t}}{\rightarrow}\real^{2n}$ with $A\times B\subseteq f(M)$.  And because it is true for all $k,t$ such that $Dm(A)<k,~Dm(B)<t$, then $Dm(A\times B)\leqslant Dm(A)+Dm(B).$
\hspace*{\fill}$\Box$ \\

\noindent
{\bf Proof of Theorem \ref{p14}.}\newline
Part 1. Let $k\in \real$ be a number such that $Dm(A)<k$, then there exists a $D^k$ function $f:B\subseteq\real^1\rightarrow \real^n$ with $A\subseteq f(B)$, and consequently there exists $M>0$ that $\forall x,y\in B \quad |f(x)-f(y)|^k<M|x-y|.$  For every $\sigma>0$ let us take a natural number $s_0$ such that $\left(\frac{M}{2^{s_0}}\right)^{1/k}<\sigma.$   Now considering the set $K^1_{s_0}$ (see (\ref{K^n_s})), we can see:
\begin{enumerate}
\item $A\subseteq\bigcup_{i\in\mathbb N}A_i$, where $A_i=f(B\cap[-i,i])$
\item $A_i\subseteq\bigcup_{\delta\in K^1_{s_0},\delta\subseteq[-i,i]}f(B\cap\delta) \quad \forall i\in\mathbb N$
\item $\forall\delta\in K^1_{s_0} \quad \diam f(B\cap\delta)\leqslant\sigma$
\item $\forall i \quad \sum_{\delta\in K^1_{s_0},\delta\subseteq[-i,i]}(\diam f(B\cap\delta))^k<2Mi~.$
\end{enumerate}

Thus for each $i \quad \dim_H A_i\leqslant k$ and therefore $\dim_H A\leqslant k$ because $\dim_H A\leqslant\sup_i \dim_H A_i$ (see 4.8 in \cite{Mattila}). And since $\dim_H A\leqslant k$ for every $k> Dm (A)$, it follows that $\dim_H A\leqslant Dm (A).$\\

\noindent
Part 2. There exist compact sets $A,B\subseteq\real^1$  such that $\dim_H A+\dim_H B<\dim_H(A\times B)$ (see \cite{Federer}), then by the Part 1 and theorem \ref{p4}  \quad $\dim_H A+\dim_H B< Dm (A)+Dm (B)$, so either $\dim_H A<Dm (A)$ or $\dim_H B<Dm (B)$ is true.
\hspace*{\fill}$\Box$ \\

\noindent
{\bf Proof of Theorem \ref{p2}.}
As follows from the Theorem \ref{p11}, it only remains to show that $Dm([0,1]^m)=m.$
In \cite{Ainouline2} it was proven that $\forall n\in \mathbb N~~ \quad \exists ~f_n:[0,1]\stackrel{D^n,{\rm onto}}{\rightarrow}[0,1]^n$~~so that $Dm([0,1]^m)\leqslant m$. Now by Theorem \ref{p14}\quad $Dm\geqslant\dim_H$, and $\dim_H[0,1]^m=m$ (see p.59 in \cite{Mattila}), then it follows that $Dm([0,1]^m)=m.$\\  
\hspace*{\fill}$\Box$ \\

\begin{lemma}\label{lemma 3.1}
Every $C^{k.\lambda}$ function $(k+\lambda\not= 0)~f:U\subseteq\real^m\rightarrow\real^n$,~with $U$~bounded, is $D^{\frac{1}{k+\lambda}}$ function on a set $\sum^*_k f=\{x\in U~;~\text{every partial derivative of}~f~\text{of order}~\leqslant k^*~\text{vanishes at}~x,~k^*=k~\text{if}~\lambda\not =0;~k^*=k-1~\text{if}~\lambda = 0\}$.
\end{lemma}

\noindent
{\bf Proof.}
If $k=0$, then conclusion of this theorem follows from the definition of $C^{0.\lambda}$ function for all $x\in U$.\newline
If $k\in\mathbb N$, we apply Taylor formulae for the points $x\in\sum^*_k f$ and use the fact that $D_{k^*}f\in C^{0.\lambda}$~~if $\lambda\not =0$~~or $D_{k^*}f\in C^1$~~if $\lambda =0$.
\hspace*{\fill}$\Box$ \\

\noindent
{\bf Proof of Theorem \ref{p1}.}
We can suppose without loss of generality that $A,B$ are bounded and $Dm(A)<Dm(B)$,  then there exist $k\in \real~:~Dm(A)<k<Dm(B)$ and a $D^k$ function $f:E\subseteq\real\stackrel{onto}{\rightarrow} A.$ By contradiction, if $A, B$ are diffeomorphic then there is a $C^1$ function $g:A\stackrel{onto}{\rightarrow} B$, and by Lemma \ref{lemma 3.1} $g\in D^1$. Now by the "Composition property" (\cite{Ainouline2}) we have a $D^k$ function $g\circ f:E\subseteq\real\stackrel{onto}{\rightarrow} B$ that contradicts $k<Dm(B).$ 
\hspace*{\fill}$\Box$ \\

\noindent
{\bf Proof of Theorem \ref{p9}.}
It follows from Main Theorem \cite{Ainouline3}, Theorem~\ref{p11}  and Lemma \ref{lemma 3.1}.
\hspace*{\fill}$\Box$\\

\noindent
{\bf Proof of Theorem \ref{p7}.}
It follows from Theorems~~\ref{p5},~\ref{p11}  and Lemma \ref{lemma 3.1}.
\hspace*{\fill}$\Box$ \\

\noindent
{\bf Proof of Theorem \ref{p15}.}
It is possible to choose a real number $p\in\real$ such that $m<p$~and ~$k<p<\frac{m}{Dm(A)}$,~then $Dm(A) < \frac{m}{p}$~ and by the Definition \ref{R-dimension} there exists a $D^{m/p}$ function $g:B\subseteq\real^1\rightarrow\real^1$~ such that $A\subseteq g(B).$\newline  
By the "$C^{<k}$-extension property"  (\cite{Ainouline3}) there exists a $C^{<[\frac{p}{m}].\frac{p}{m}-[\frac{p}{m}]}$-function ~$q:\real^1\longrightarrow\real^1$ such that $A\subseteq q(\sum_{p/m}q)$,~((Designation: We write {\boldmath $F\in C^{<k.\lambda}$} if $F\in C^{p.\beta}$~ for all $p+\beta<k+\lambda$))  then  by the Main Theorem \cite{Ainouline2} there exists a $C^{<[p].p-[p]}$-function $f:\real^m\rightarrow\real^1$ such that $q(\sum_{p/m}q)$ is the critical values set of the function $f$(Note that $C^{k+\lambda}\subseteq C^{k.\lambda}$).
From $A\subseteq q(\sum_{p/m}q)$ and $k<p$, we conclude that $f\in C^{[k].k-[k]}$~ and $A$ is subset of the critical values set of $f$.
\hspace*{\fill}$\Box$\\

\begin{lemma}\label{g(B)}
$\forall k>1$~~there exists a $D^k$-function $g:B\subseteq\real^1\stackrel{{onto}}{\longrightarrow}\real^1$~such that \meas(B)=0.
\end{lemma}

\noindent
{\bf Proof.}
We define a family of sets $P=\bigcup_{s\in\mathbb N}P_s$ on $\real^1.$  For $s=1 \quad P_1=\{G_{i_0,i_1}~:~i_0\in\mathbb N,~i_1=1,2\}$, where each $G_{i_0,i_1}$  is a closed interval of length $\frac{1}{2^k}$~ located in the middle of $Q^1_{i_0,i_1}\in K^1_1$ defined in 
(\ref{K^n_s}). And if $P_{s-1}=\{G_{i_0,\dots,i_j,\dots,i_{s-1}}:i_j=1,2\}$ have already been constructed, then $P_s=\{G_{i_0,\dots,i_j,\dots,i_s}~:~i_j=1,2\}$,~ where $G_{i_0,\dots,i_{s-1},1},G_{i_0,\dots,i_{s-1},2}$ are two closed intervals of length $\frac{1}{2^{ks}}$, first located at the beginning, second at the end of $G_{i_0,\dots,i_{s-1}}$. We can notice that distance between them is  
\begin{equation}\label{distance}
\frac{1}{2^{k(s-1)}}-\frac{2}{2^{ks}}=(2^k-2)\frac{1}{2^{ks}}.
\end{equation}

Let $B\subseteq \real^1$  be the Cantor set defined by the family of sets $P$. Let us define a function $g:B\subseteq\real^1\stackrel{onto}{\longrightarrow}\real^1$ with the condition
$g(G_{i_0,\dots,i_s}\cap B)=Q^1_{i_0,\dots,i_s}\in K^1_s$ ~(see (\ref{K^n_s})).  To check that $g\in D^k$, let us take arbitrary points $x,y\in B,~~(x<y)$, then there exists $s\in\mathbb N$ such that $x\in G_{i_0,\dots,i_{s-1},1},~y\in G_{i_0,\dots,i_{s-1},2}$.

\begin{equation}\label{property of g}
|g(x)-g(y)|^k\leqslant\left(\frac{1}{2^{s-1}}\right)^k\leqslant K|x-y|
\end{equation}
where $K=\frac{2^k}{2^k-2}$ depends only on k.

To prove that $\meas(B)=0$, suffice it to show that $\meas(B\cap[a,a+1])=0 \quad\forall a\in\mathbb Z$. And the $\meas(B\cap[a,a+1])\leqslant \lim\limits_{s\to\infty}2^s\cdot\frac{1}{2^{ks}}=0$ \quad  because $k>1$.
\hspace*{\fill}$\Box$\\

\noindent
{\bf Proof of Theorem \ref{p21}.}
Let $a,b\in\real$ be such that $Dm(A)< b < a$, then by the definition \ref{R-dimension} there exists a $D^b$-function $f_a:F\subseteq\real^1\rightarrow\real^n$~~ such that $A\subseteq f_a(F)$. We can suppose that the coefficient $K$ required for $f_a$ by the definition \ref{D^k}  is equal to 1 (stretching $F$ if necessary). By Lemma \ref{g(B)} and because $\frac{a}{b}>1$, there exists $g:G\subseteq\real^1\stackrel{onto}{\longrightarrow}\real^1$  such that\quad $\meas(G)=0,\quad\exists K^*>0~~:~~\forall x,y\in G \quad |g(x)-g(y)|^{a/b}\leqslant K^*|x-y|.$

Let us define a function $q:Q\subseteq\real^1\rightarrow\real^1$  such that
$Q=\{K^*x~;~x\in G\} \quad\forall z\in Q \quad q(z)=g(\frac{z}{K^*})$, and a function $h=f_a\circ q\upharpoonright H$, where $H=Q\cap q^{-1}(F).$

Now we have the following:
\begin{enumerate}
\item $A\subseteq h(H)$
\item $\meas(H)=0$
\item $\forall x,y\in H\quad |h(x)-h(y)|^a~=~\left|f_a(q(x))-f_a(q(y))\right|^{ab/b}~\leqslant~\left|q(x)-q(y)\right|^{a/b}~=\newline
 \quad~=~\left|g(\frac{x}{K^*})-g(\frac{y}{K^*})\right|^{a/b}~\leqslant~K^*\left|\frac{x}{K^*}-\frac{y}{K^*}\right|~=~1\cdot|x-y|.$
\end{enumerate}

And by the definition \ref{m_a} and from 1,2,3 it follows that $m_a(A)=0.$
\hspace*{\fill}$\Box$\\

\begin{lemma}\label{lem3.2*}
Let $A=\bigcup_{i\in\mathbb N}A_i\subseteq\real^n,~n\in\mathbb N,~k\geqslant 0,~k\in\real$ and $\{f_i:B_i\subseteq[a_i,b_i]\rightarrow\real^n,~i\in\mathbb N\}$ be such that ~$\forall i\in \mathbb N  \quad  A_i\subseteq f_i(B_i)$, and $\forall b,b^\prime\in B_i \quad|f_i(b)-f_i(b^\prime)|^k\leqslant |b-b^\prime|$. Then there exists $f:B\subseteq \real^1\rightarrow\real^n$ with $A\subseteq f(B)$,such that $\forall b,b^\prime\in B \quad|f(b)-f(b^\prime)|^k\leqslant |b-b^\prime|$ ~and $\meas(B)\leqslant \sum_{i\in \mathbb N}\meas(B_i)$.
\end{lemma}

\noindent
{\bf Proof:} Similar to the prove of lemma \ref{lemma 3.2}.
\hspace*{\fill}$\Box$\\

\noindent
{\bf Proof of Theorem \ref{p22}.}
1)and 2)~~are obvious.\quad
3)~~Let's consider that $\forall i\in \mathbb N \quad m_a(A_i)<+\infty$~~~(otherwise it is obvious).
Now for an arbitrary fixed $\epsilon >0$ take a set $\{\epsilon_i>0:i\in \mathbb N \}$~~such that $\sum_{i\in \mathbb N}\epsilon_i \leqslant \epsilon$.  Then by the definition  ~ \ref{m_a}~~there exists a family of function $\{F_i:B_i\subseteq \real^1\rightarrow\real^n,~i\in\mathbb N\}$~~~  such that~~~ $\forall i\in \mathbb N  \quad  A_i\subseteq F_i(B_i)$, \quad  $\forall b,b^\prime\in B_i \quad |F_i(b)-F_i(b^\prime)|^a\leqslant |b-b^\prime|$,~~ and~~ $meas (B_i)\leqslant m_a(A_i)+\epsilon_i$.
Let's consider a denumerable family of function $\{f_{ij}=F_i\upharpoonright B_{ij} ; B_{ij}=B_i\cap [j,j+1],i\in \mathbb N,j\in \mathbb Z\}$. Then from the lemma \ref{lem3.2*} there exists $f:B\subseteq \real^1\rightarrow\real^n$ with $\bigcup_{i\in \mathbb N}A_i=\bigcup_{i\in \mathbb N}\bigcup_{j\in \mathbb Z}f_{ij}(B_{ij})\subseteq f(B)$  such that $\forall b,b^\prime\in B \quad |f(b)-f(b^\prime)|^a\leqslant |b-b^\prime|$ ~and $\meas(B)\leqslant \sum_{i\in \mathbb N}\sum_{j\in \mathbb Z}B_{ij}\leqslant \sum_{i\in \mathbb N}\meas(B_i)\leqslant \sum_{i\in \mathbb N}m_a(A_i)~+~\epsilon$. And because it is true for any $\epsilon >0$, by the definition \ref{m_a} we have that $m_a(\bigcup_{i\in \mathbb N}A_i)\leqslant\sum_{i\in \mathbb N}A_i$.
\hspace*{\fill}$\Box$\\

\end{document}